\documentclass[reqno,10pt]{article}
\usepackage{amssymb}
\usepackage{color}

\begin{document}

\newtheorem{lem}{Lemma}

\newenvironment{pf}{\noindent {\it Proof.}\,\,}{{\hfill $\square$} \vskip 6pt}

\title{Representing Primes as $x^2 + 5y^2$: An Inductive Proof that Euler Missed} %
\author{Ying Zhang} %

\date{October 3, 2006}

\maketitle

\section{HISTORICAL INTRODUCTION}\label{s:1}

In this note we present an elementary inductive proof which Euler
could have obtained, for his assertion that every prime of the form
$20n+1$ or $20n+9$ is a sum $x^2 + 5y^2$, had he refined a bit his
proof for Fermat's theorem that every prime of the form $4n+1$ is a
sum of two squares.

Here and throughout this note all letters are assumed to stand for
{\it nonnegative} integers, unless otherwise specified.
It is our pleasure to start by briefly reviewing the story told by
Cox in the nice book \cite{cox1989book}.

\vskip 3pt

Pierre de Fermat (1601--1665), who had done pioneering work on
representing primes as $x^2 + ny^2$, stated, but did not write down
a proof, that he had proved by his favorite method of infinite
descent the following:

\vskip 3pt

(i)  Every prime  $p \equiv 1 \,\,{\rm mod}\,\, 4$ is a sum $x^2 +
y^2$.

\vskip 2pt

(ii)  Every prime  $p \equiv 1, 3 \,\,{\rm mod}\,\, 8$ is a sum $x^2
+ 2y^2$.

\vskip 2pt

(iii)  Every prime $p \equiv 1 \,\,{\rm mod}\,\, 3$  is a sum $x^2 +
3y^2$.

\vskip 3pt

\noindent He also conjectured but could not prove that

\vskip 3pt

(iv)  The product of two primes, each of which is $\equiv 3, 7
\,\,{\rm mod}\,\, 20$, is a sum $x^2 + 5y^2$.

\vskip 3pt

Leonhard Euler (1707--1783) heard of Fermat's results and spent 40
years proving (i)--(iii) and considering their generalizations. This
finally led him to the discovery of the quadratic reciprocity,
although he could not provide a solid proof for it. By working out
numerous examples on representing primes as $x^2+ny^2$ for various
$n$, he discovered more patterns. Some of his discoveries which he
could not prove are:

\vskip 3pt

(v)  Every prime  $p \equiv 1, 9 \,\,{\rm mod}\,\, 20$ is a sum $x^2
+ 5y^2$.

\vskip 2pt

(vi)  For every prime  $p \equiv 3, 7 \,\,{\rm mod}\,\, 20$, $2p$ is
a sum $x^2 + 5y^2$.

\vskip 2pt

(vii)  A prime  $p=x^2 + 27y^2 $ if and only if $ p \equiv 1
\,\,{\rm mod}\,\, 3$ and 2 is a cubic residue modulo $p$.

\vskip 2pt

(viii)  A prime  $p=x^2 + 64y^2 $ if and only if $ p \equiv 1
\,\,{\rm mod}\,\, 4$ and 2 is a biquadratic residue modulo $p$.

\vskip 3pt

Joseph-Louis Lagrange (1736--1813) and Adrien-Marie Legendre
(1752--1833) later developed the form theory as well as the genus
theory to prove (iv)--(vi). Indeed, they could prove (v) and that

\vskip 3pt

(v$'$)  Every prime  $p \equiv 3, 7 \,\,{\rm mod}\,\, 20$ is a sum
$2x^2 + 2xy + 3y^2$ (where one of $x,\, y$ may be negative).

\vskip 3pt

\noindent Then (iv) and (vi) follow immediately from the following
two identities:
\begin{eqnarray*}
(2x^2 + 2xy +3y^2)(2a^2 + 2ab + 3b^2)\!\!\!\!\!&=&\!\!\!\! %
(2ax + bx +ay +3by)^2 + 5(bx - ay)^2; \\ %
2(2x^2 + 2xy +3y^2)\!\!\!\!\!&=&\!\!\!\!(2x + y)^2 + 5y^2.
\end{eqnarray*}
Both Legendre and Lagrange, however, could prove neither (vii) nor
(viii).

\vskip 3pt

It was Carl Friedrich Gauss (1777--1855) who had finally tackled
(vii) as well as (viii) using his cubic and biquadratic
reciprocities. And, years before this, it was also Gauss who gave
the first rigorous proof of the quadratic reciprocity. The
interested reader is referred to \cite{cox1989book} to enjoy the
rest of the story.

What we shall show is that in fact Euler {\it could} have proved
(iv)--(vi) had he just refined his proof of (i)--(iii), and hence
the above told story would be somewhat different.

\section{A REVIEW OF EULER'S PROOF}\label{s:2}

Let us first briefly review Euler's proof.
According to [1], a version of Euler's proof of (i) goes as follows.
For prime $p \equiv 1 \,\,{\rm mod}\,\, 4$, there is an $r < p$ such
that $pr = x^2 + 1$. For each prime factor $q$ of $r$, since $-1$ is
a quadratic residue modulo $q$, it follows that either $q \equiv 1
\,\,{\rm mod}\,\, 4 $ or $q=2$. We assume by induction that each
such $q$ is a sum of two squares. Then a cancelation lemma, Lemma
\ref{lem:1} below, enables to cancel the prime factors of $r$ one by
one. As a result, one obtains a representation of $p$ as a sum of
two squares.

By a $(1, n)$-{\it representation} we mean an expression of the from
$x^2 + ny^2$. The following lemma appears as Lemma 1.4 in
\cite{cox1989book}. Here, for the convenience of the reader, we
include the proof given in \cite{cox1989book} in a slightly
shortened form.

\begin{lem}\label{lem:1}
Suppose $p$ and $pr$ each has a $(1,n)$-representation, where $p$ is
a prime. Then $r$ has a $(1,n)$-representation.

\end{lem}

\begin{pf} Suppose $p = a^2 + nb^2$, and $pr = x^2 + ny^2$.
Then
\begin{eqnarray*}
p^{2}r = (ax \pm nby)^2 + n(ay \mp bx)^2.
\end{eqnarray*}
Note that $p \,\mid \,(ay
- bx)(ay + bx)$ since
\begin{eqnarray*}
(ay - bx)(ay + bx) = (a^2 + nb^2)y^2 - b^2(x^2 + ny^2).
\end{eqnarray*}
It follows that either $p \mid ay - bx$ or $p \mid ay + bx$;
correspondingly, $p \mid ax + nby$ or $p \mid ax - nby$.
Consequently, we have one of the following holds:
\begin{eqnarray*}
&& r = ((ax + nby)/p)^2 + n ((ay - bx)/p)^2, \\ %
&& r = ((ax - nby)/p)^2 + n ((ay + bx)/p)^2. %
\end{eqnarray*}
This proves Lemma \ref{lem:1}.
\end{pf}

It is not hard to see that Euler's proof also applies to the cases
of $x^2 + 2y^2$ and $x^2 + 3y^2$ after minor modifications. This is
because in the representation $pr = x^2 + 2$ (resp. $pr = x^2 + 3$)
all of the prime factors of $r$ (with one or two exceptions which
are easy to deal with) are of the desired type, so we can again use
Lemma \ref{lem:1} and the inductive hypothesis to cancel them one by
one.

\vskip 1pt

However, it is {\it not} the case for $x^2 + 5y^2$. Note that for
$pr = x^2 + 5$, where we may assume that $5 \nmid r$, each prime
factor $q$ of $r$ is such that either $q \equiv 1, 9, 3, 7 \,\,{\rm
mod}\,\, 20$ or $q=2$, hence  Lemma \ref{lem:1} is not enough to
cancel all prime factors $q$ of $r$; we have to deal with those $q
\equiv 3, 7 \,\,{\rm mod}\,\, 20$ and $q=2$. Indeed we {\it do} have
such a cancelation lemma (Lemma \ref{lem:2} in \S \ref{s:3}) which
enables us, with the help of a small trick, to cancel such factors
pair by pair under the inductive hypothesis that (iv)--(vi) hold for
all primes $q < p$ such that $q \equiv 1, 9, 3, 7 \,\,{\rm mod}\,\,
20$ or $q = 2$. It turns out that we must prove (iv)--(vi) {\it
simultaneously} by induction.
The rest of this note consists of the detailed statements and
proofs.

\section{A SECOND CANCELATION LEMMA}\label{s:3}

For the convenience of further exposition, we make the following
definition.

\vskip 6pt

\noindent {\bf Definition.}\, A $(1, n)$-representation is said to
be {\it nontrivial} if both of $x$ and $y$ are nonzero; it is {\it
proper} if $x$ and $y$ are relatively prime.

\vskip 6pt

\noindent {\bf Remarks.}\, The following items ({\bf a})--({\bf f})
can be easily checked.

\vskip 2pt

({\bf a}) A proper $(1, n)$-representation is automatically
nontrivial unless it equals $1$ or $n$.

\vskip 2pt

({\bf b}) A $(1, n)$-representation of a prime $p$, where $p \neq
n$, is automatically proper and nontrivial.

\vskip 2pt

({\bf c}) A nontrivial $(1, n)$-representation of the product of two
primes is always proper.

\vskip 2pt

({\bf d}) If $p$ is a prime such that $p \nmid r$ and $p \nmid  n$,
then any $(1, n)$-representation of $pr$ is nontrivial.

\vskip 2pt

({\bf e}) There is the following very useful {\it Euler identity}
which expresses the product of two $(1, n)$-representations as a
$(1, n)$-representation in two ways:
\begin{eqnarray}\label{eqn:eulerid}
(a^2 + nb^2)(x^2 + ny^2) = (ax \pm nby)^2 + n(ay \mp bx)^2.
\end{eqnarray}

\vskip 2pt

({\bf f}) If an odd $s$, where $5 \nmid s$, has a nontrivial, proper
$(1,5)$-representation $s=a^2 + 5b^2$, then
\begin{eqnarray*}
s^2 = (a^2 - 5b^2)^2 + 5(2ab)^2
\end{eqnarray*}
is a nontrivial, proper $(1,5)$-representation, since every prime
common factor of $a^2 - 5b^2$ and $2ab$ is a common factor of $a$
and $b$.

\vskip 6pt

\noindent Our second cancelation lemma, Lemma \ref{lem:2} below,
enables us to cancel $q^2$, where $q$ is an odd prime and $q^2$ has
a nontrivial $(1, n)$-representation, from a given $(1,
n)$-representation of $q^2r$ and obtain a $(1, n)$-representation of
$r$.

\begin{lem}\label{lem:2}
Suppose $q^2$ has a nontrivial $(1, n)$-representation, where $q$ is
an odd prime. If $q^2 r$ has a $(1, n)$-representation, then $r$ has
a $(1, n)$-representation. Moreover, if $q^2 r$ has a nontrivial,
proper $(1, n)$-representation and $r \neq 1,n$, then $r$ has a
nontrivial, proper $(1, n)$-representation.
\end{lem}

\begin{pf} Suppose $q^2 = a^2 + nb^2$ is
nontrivial. Then it is proper by Remark ({\bf c}). Let $q^2 r = x^2
+ ny^2$. From the Euler identity (\ref{eqn:eulerid}), we have
\begin{eqnarray*}
q^4 r = (ax \pm nby)^2 + n(ay \mp bx)^2.
\end{eqnarray*}
Then $q^2 \mid (ay - bx)(ay + bx)$ as in the proof of Lemma
\ref{lem:1}.

\vskip 1pt

First, suppose $q^2 \nmid ay - bx$ and $q^2 \nmid ay + bx $. Then we
must have $q \mid ay - bx$ and $q \mid ay + bx$. Hence $q \mid 2ay$.
Since $q$ is odd, $q \mid ay$. We then have $q \mid y$ since $q \,
\nmid \, a$ (otherwise $a = 0$, or $a = q$ and $b = 0$, a
contradiction). Similarly, $q \mid x$. Consequently, we have
$$
r = (x/q)^2 + n(y/q)^2.
$$
In this case $q^2 r = x^2 + ny^2$ is not proper.

Now we may suppose \,$q^2 \mid ay - bx$ \,or\, $q^2 \mid ay + bx$.
Then \,$q^2 \mid ax + nby$ \,or\, $q^2 \mid ax - nby$ \,accordingly.
Consequently, we have one of the following two holds:
\begin{eqnarray}
&& r = ((ax + nby)/{q^2})^2 + n ((ay - bx)/{q^2})^2, \label{eqn:case1} \\%
&& r = ((ax - nby)/{q^2})^2 + n ((ay + bx)/{q^2})^2. \label{eqn:case2} %
\end{eqnarray}

\noindent {\sc Claim I.} \ The above obtained $(1,
n)$-representation of $r$ is nontrivial and proper if so is $q^2 r =
x^2 + ny^2$.

\vskip 3pt

\noindent {\sc Proof of Claim I.} \ In the case where $q^2 \mid ay -
bx$ and $q^2 \mid ax + nby$, it follows from the identities
\begin{eqnarray*}
&&x = a ((ax + nby)/{q^2}) - nb\, ((ay - bx)/{q^2}), \\%
&&y = b\, ((ax + nby)/{q^2}) + a ((ay - bx)/{q^2}) %
\end{eqnarray*}
that $(ax + nby)/{q^2}$ and $(ay - bx)/{q^2}$ are relatively prime
since so are $x$ and $y$. Hence the $(1, n)$-representation
(\ref{eqn:case1}) is proper.

In the case where $q^2 \mid ay + bx$ and $q^2 \mid ax - nby$, it
follows from the identities
\begin{eqnarray*}
&&x = a ((ax - nby)/{q^2}) + nb\, ((ay + bx)/{q^2}), \\%
&&y = -b\, ((ax - nby)/{q^2}) + a ((ay + bx)/{q^2}) %
\end{eqnarray*}
that $(ax - nby)/{q^2}$ and $(ay + bx)/{q^2}$ are relatively prime
since so are $x$ and $y$. Hence the $(1, n)$-representation
(\ref{eqn:case2}) is proper.

Since $r \neq 1,n$, in either case the $(1, n)$-representation of
$r$ is nontrivial by Remark ({\bf a}). Claim I is thus proved.

\vskip 3pt

This completes the proof of Lemma \ref{lem:2}.
\end{pf}

\noindent The case where $q=2$ and $n=5$ is simple and is considered
in the following

\vskip 6pt

\noindent {\bf Addendum to Lemma \ref{lem:2}.}\, If $2^2r$ has a
$(1, 5)$-representation $2^2r=x^2 + 5y^2$, then $x$ and $y$ must be
both even (by a simple modulo $4$ argument), and hence $r = (x/2)^2
+ 5(y/2)^2$. Furthermore, if $x^2 + 5y^2$ is nontrivial, so is
$(x/2)^2 + 5(y/2)^2$.

\section{THE PROOF THAT EULER MISSED}\label{s:4}

For convenience of later reference, we restate the assertions
(iv)--(vi) in \S \ref{s:1} as

\vskip 6pt

\noindent {\bf Theorem.} \,\, {\it {\rm(1)} Every prime $p \equiv 1,
9 \,\,{\rm mod}\,\, 20$ has a $(1, 5)$-representation.

\vskip 1pt

{\rm(2)} For every pair of primes $q,q'$ such that $q \equiv 3,7
\,\,{\rm mod}\,\, 20$ and either $q'\equiv 3,7 \,\,{\rm mod}\,\, 20$
or $q'=2$, their product $q q'$ has a nontrivial $(1,
5)$-representation.}

\vskip 6pt

\noindent It is the following inductive proof that Euler missed.

\vskip 6pt

\begin{pf} Suppose by induction that (1) and (2) hold for all primes
$p, q, q'$ which are less than a certain prime $\pi$ where $\pi
\equiv 1,9,3,7 \,\,{\rm mod}\,\, 20$. We need to show that

\vskip 3pt

$(1)_{\pi}$ \ If $\pi \equiv 1,9 \,\,{\rm mod}\,\, 20$ then $\pi$
has a $(1, 5)$-representation.

$(2)_{\pi}$ \ If $\pi \equiv 3,7 \,\,{\rm mod}\,\, 20$ then for
every prime $q \le \pi$ such that $q \equiv 3,7$ $\,\,{\rm mod}\,\,
20$ or $q = 2$, $\pi q$ has a nontrivial (hence proper) $(1,
5)$-representation.

\vskip 3pt

To start, we have from the quadratic reciprocity that for a prime $p \neq 2, 5$, %
\vskip 3pt %
\centerline{$-5$ is a quadratic residue mod $p$
\quad $\Longleftrightarrow$ \quad $p \equiv 1,9,3,7 \,\,{\rm mod}\,\, 20$.} %
\vskip 3pt %
\noindent Hence there is a $(1,5)$-representation
\begin{eqnarray}\label{eqn:pir=}
\pi r =\, x^2 + 5y^2.
\end{eqnarray}
Here (\ref{eqn:pir=}) initially holds for some $x \le (\pi - 1)/2$
and $y=1$; it follows that $x^2 + 5y^2 < {\pi}^2$ and hence $r <
\pi$. After reduction if necessary, it can be assumed that $5 \nmid
r$ and that (\ref{eqn:pir=}) is a nontrivial, proper
$(1,5)$-representation.

\vskip 4pt

\noindent {\sc Claim II.}\, When $r$ in (\ref{eqn:pir=}) is
minimized, we have either $r = 1$ or $r = q'$, where $q'$ is a prime
such that either $q' \equiv 3,7 \,\,{\rm mod}\,\, 20$ or $q'= 2$.

\vskip 4pt

\noindent {\sc Proof of Claim II.}\, Since $-5$ is a quadratic
residue mod $p$, for each prime factor $q$ of $r$, we have $q< \pi$
and either $q\equiv 1,9,3,7 \,\,{\rm mod}\,\, 20$ or $q= 2$. Our
idea is to manage to cancel the prime factors of $r$ one by one for
those congruent to $1,9$ modulo $20$, and pair by pair for those
congruent to $3,7$ modulo $20$ or equal to $2$.

If $r$ has a prime factor $q$ such that $q \equiv 1,9 \,\,{\rm
mod}\,\, 20$, then, by the inductive hypothesis, $q$ has a $(1,
5)$-representation. By Lemma \ref{lem:1}, $\pi r'$, where $r' =
r/q$, has a $(1, 5)$-representation. Hence, by minimizing $r$ in
(\ref{eqn:pir=}), we may assume that $r$ has no prime factors
congruent to $1,9$ modulo $20$.

Now each prime factor $q$ of $r$ is of the form %
either $q \equiv 3,7 \,\,{\rm mod}\,\, 20$ or $q = 2$. %
If the number of such prime factors of $r$, counted  with
multiplicity, is at least $2$, let $q, q'$ be two of them and set
$r' = r/(q q')$.

If $q = q'$ then $qq'$ has a nontrivial $(1, 5)$-representation by
the inductive hypothesis, hence we can apply Lemma \ref{lem:2} and
its addendum directly to cancel $q^2$ from $\pi r=(\pi r')(q^2)$ and
obtain a $(1, 5)$-representation of $\pi r'$.

If $q \neq q'$ then $q q'$ has a nontrivial $(1, 5)$-representation
by the inductive hypothesis again. Now $q^2 (q')^2 \pi r' = (q
q')(\pi r)$ has a $(1, 5)$-representation by the Euler identity
(\ref{eqn:eulerid}), and applying Lemma \ref{lem:2} and its addendum
{\it twice} implies that $\pi r'$ has a $(1,
5)$-representation---here is the trick used.
This finishes the proof of Claim II.

\vskip 3pt

We proceed to prove the inductive step. By minimizing $r$ in
(\ref{eqn:pir=}), we are in one of the alternatives described in
Claim II.

First, we prove {$(1)_{\pi}$}. In this case $\pi \equiv 1,9 \,\,{\rm
mod}\,\, 20$. One must have $r = 1$ and hence $\pi$ has a $(1,
5)$-representation; otherwise, $r = q' \equiv 3,7 \,\,{\rm mod}\,\,
20$ or $r=q'=2$, but then $\pi q' = x^2 + 5y^2 \equiv 3,7,\pm 2
\,\,{\rm mod}\,\, 20$ and consequently $x^2 \equiv \pm 2 \,\,{\rm
mod}\,\, 5$, a contradiction. This proves {$(1)_{\pi}$}.

To prove {$(2)_{\pi}$}, suppose $\pi \equiv 3,7 \,\,{\rm mod}\,\,
20$. One must have $r = q'$ as described in Claim II; otherwise $r =
1$, which implies that $\pi r = x^2 + 5y^2 \equiv 3,7 \,\,{\rm
mod}\,\, 20$, a contradiction.
Thus $\pi q'$ has a $(1, 5)$-representation, which is automatically
nontrivial and proper, for some prime $q' < \pi$ such that either
$q'\equiv 3,7 \,\,{\rm mod}\,\, 20$ or $q'=2$. Then Remark ({\bf f})
implies that ${\pi}^2 (q')^2 = (\pi q')(\pi q')$ has a nontrivial,
proper $(1, 5)$-representation. By the inductive hypothesis, either
$q' = 2$ or $(q')^2$ has a nontrivial $(1, 5)$-representation. Lemma
\ref{lem:2} and its addendum then give a nontrivial, proper $(1,
5)$-representation of ${\pi}^2$.

To prove the remaining part of {$(2)_{\pi}$}, let $q < \pi$ be any
prime such that either $q \equiv 3,7 \,\,{\rm mod}\,\, 20$ or $q=2$.
Then either $q = q' = 2$ or, by the inductive hypothesis, $q q'$ has
a nontrivial $(1, 5)$-representation. Thus $\pi q(q')^2 = (\pi
q')(q' q)$ has a $(1, 5)$-representation by the Euler identity
(\ref{eqn:eulerid}). On the other hand, by the inductive hypothesis,
either $q' = 2$ or $(q')^2$ has a nontrivial $(1,
5)$-representation. Now Lemma \ref{lem:2} and its addendum give a
$(1, 5)$-representation for $\pi q$, which is automatically
nontrivial. This proves {$(2)_{\pi}$}.

The theorem is thus proved by induction.
\end{pf}


\noindent Note that we have proved (iv)--(vi) without reference to
(v$'$). More interesting is that in fact (v$'$) follows from (vi).
To see this, for any prime $p \equiv 3,7 \,\,{\rm mod}\,\, 20$, let
$2p = x^2 + 5y^2$. It follows that both $x$ and $y$ are odd. Hence
$x = 2x' + y$, where $x'$ may be negative. Then $p = 2{x'}^2 + 2x'y
+3y^2$ gives a desired representation.


Among many other existing elementary proofs of Fermat's theorem (i),
we cannot help but mention Zagier's beautiful ``one-sentence
proof\,'' (see \cite{zagier1990amm}, or as explained in
\cite{aigner-ziegler2004book}) to conclude this note.

\vskip 10pt

\noindent {\bf ACKNOWLEDGEMENTS.} \,\, The author would like to
thank H. Y. Loke for teaching him Number Theory years ago and for
encouragement. Thanks are also due to the referees whose
constructive suggestions helped improve the exposition of this note.
The author is supported by a CNPq-TWAS postdoctoral fellowship and
partially by NSFC grant No.10671171.

\vskip 10pt

\noindent {Department of Mathematics, Yangzhou University, Yangzhou 225002, CHINA} %

\noindent E-mail: \verb"yingzhang@yzu.edu.cn"

\vskip 6pt

\noindent {{\sc Current Address:} IMPA, Estrada Dona Castorina 110, Rio de Janeiro 22460, BRAZIL %
\quad E-mail: \verb"yiing@impa.br"}

\end{document}